\RequirePackage{ifpdf}
\ifpdf 
\documentclass[pdftex]{sigma}
\else
\documentclass{sigma}
\fi

\begin{document}
\allowdisplaybreaks

\renewcommand{\thefootnote}{$\star$}

\renewcommand{\PaperNumber}{037}

\FirstPageHeading

\ShortArticleName{Hilbert Transforms Associated with Dunkl--Hermite Polynomials}

\ArticleName{Hilbert Transforms Associated\\ with Dunkl--Hermite Polynomials\footnote{This paper is a contribution to the Special
Issue on Dunkl Operators and Related Topics. The full collection
is available at
\href{http://www.emis.de/journals/SIGMA/Dunkl_operators.html}{http://www.emis.de/journals/SIGMA/Dunkl\_{}operators.html}}}

\Author{N\'ejib BEN SALEM and Taha SAMAALI}

\AuthorNameForHeading{N.~Ben Salem and T.~Samaali}

\Address{Department of Mathematics, Faculty of Sciences of Tunis, \\
Campus Universitaire, 2092  Tunis, Tunisia}
\Email{\href{mailto:Nejib.BenSalem@fst.rnu.tn}{Nejib.BenSalem@fst.rnu.tn}}

\ArticleDates{Received October 14, 2008, in f\/inal form March 12,
2009; Published online March 25, 2009}

\Abstract{We consider expansions of functions in $L^{p}(\mathbb{R},|x|^{2k}dx)$, $1\leq p<+\infty$ with respect to Dunkl--Hermite functions in the rank-one setting.
We actually def\/ine the heat-dif\/fusion and Poisson integrals in the one-dimensional Dunkl setting and study their properties. Next, we def\/ine and deal with Hilbert transforms and conjugate Poisson integrals in the same setting. The formers occur to be Calder\'{o}n--Zygmund operators and hence their mapping properties follow from general results.}


\Classification{42A50; 42C10}

\renewcommand{\thefootnote}{\arabic{footnote}}
\setcounter{footnote}{0}

\section{Introduction}
The study of Dunkl operators has known a considerable growth during the last two decades due to their relevance in mathematical physics and since they give the way to built a parallel to the theory of spherical Harmonics based on f\/inite root systems and depending on a set of real parameters. In this spirit, the Hilbert transform, a basic tool in signal processing and in Fourier and harmonic analysis as well, may be def\/ined by means of partial derivatives, so that, since the commutative algebra of Dunkl operators generalize the one of partial derivatives, it is natural to extend the study of Hilbert transforms and connected topics as heat dif\/fusion, Poisson integrals and others to the Dunkl setting.
In this work, we start with investigating the rank-one case, that is why we sketch some facts that are subsequently needed. Let $k$ be a nonnegative parameter and let $T_{k}$ be the Dunkl operator acting on smooth functions $f$ as
\[
T_{k}f(x)=f'(x)+k\dfrac{f(x)-f(-x)}{x}, \qquad f\in C^{1}({\mathbb R}).
\]
To that operator is associated the so-called Dunkl--Hermite operator on ${\mathbb R}$  denoted $L_{k}$ and def\/ined by
\[
L_{k}=T_{k}^{2}-x^{2}.
\]
Its spectral decomposition is given by the Dunkl--Hermite functions $h_{n}^{k}$  def\/ined by
\[
h_{n}^{k}(x)=e^{-\frac{x^{2}}{2}}H_{n}^{k}(x),
\]
 where  $H_{n}^{k}$  are the generalized Hermite polynomials which we call the Dunkl--Hermite polynomials as in~\cite{ElGar}, namely (see~\cite{Rosl0})
\[
L_{k}h_{n}^{k}(x)=-(2n+2k+1)h_{n}^{k}(x).
\]
 Recall also that $H_n^k$ were introduced in  \cite{Rosen} and studied by R\"{o}sler in \cite{Rosl0}, whence
\begin{gather*}
H_{2n}^{k}(x)=(-1)^{n}\sqrt{\frac{n!}{\Gamma(n+k+\frac{1}{2})}}L_{n}^{k-\frac{1}{2}}\left(x^{2}\right),
\\
H_{2n+1}^{k}(x)=(-1)^{n}\sqrt{\frac{n!}{\Gamma(n+k+\frac{3}{2})}}xL_{n}^{k+\frac{1}{2}}\left(x^{2}\right),
\end{gather*}
where   $L_{n}^{\alpha}$  are the Laguerre polynomials of index  $\alpha\geq -\frac{1}{2}$, given by
\[
L_{n}^{\alpha}(x)=\frac{1}{n!}x^{-\alpha}e^{x}\frac{d^{n}}{dx^{n}}\left(x^{n+\alpha}e^{-x}\right).
\]
It is well known that the system $\{{H_{n}^{k}}\}_{n\geq0}$ is complete and orthonormal  in $L^{2}({\mathbb R},e^{-x^{2}}|x|^{2k}dx)$, therefore the system $\{h_{n}^{k}\}_{n\geq0}$ is complete and orthonormal in $L^{2}({\mathbb R},|x|^{2k}dx)$.

  Hereafter, $L^{p}({\mathbb R},|x|^{2k}dx), 1\leq p<+\infty$ is the space of measurable functions on ${\mathbb R}$ satisfying
\[
\|f\|_{k,p}:=\left(\int_{{\mathbb R}}|f(x)|^{p}|x|^{2k}dx\right)^{\frac{1}{p}}<+\infty,
\]
and $f$  belongs to  $L^{p}({\mathbb R},|x|^{2k}dx)$, $1\leq p<+\infty,$ unless mentioned. For a given $f$, one def\/ines the heat-dif\/fusion integral $G_{k}(f)$  by
\[
G_{k}(f)(t,x)=\sum_{n=0}^{+\infty}e^{-t(2n+2k+1)}a_{n}^{k}(f)h_{n}^{k}(x), \qquad t > 0 , \qquad x\in{\mathbb R},
\] where
\[
a_{n}^{k}(f)=\int_{{\mathbb R}}f(t)h_{n}^{k}(t)|t|^{2k}dt.
\]
We shall establish that  $G_{k}(f)$  possesses the following integral representation
\[
G_{k}(f)(t,x)=\int_{{\mathbb R}}P_{k}(t,x,y)f(y)|y|^{2k}dy,
\]
where
\[
P_{k}(t,x,y)=\sum_{n=0}^{+\infty}e^{-t(2n+2k+1)}h_{n}^{k}(x)h_{n}^{k}(y).
\]
We shall prove that $G_{k}(f)(t,\cdot)$, $t>0$, satisf\/ies
\[
\|G_{k}(f)(t,\cdot)\|_{k,p}\leq (\cosh(2t))^{-(k+\frac{1}{2})}\|f\|_{k,p}.
\]
Next, we def\/ine the Poisson integral $F_{k}(f)$  by
\[
F_{k}(f)(t,x)=\sum_{n=0}^{+\infty}e^{-t\sqrt{2n+2k+1}}a_{n}^{k}(f)h_{n}^{k}(x), \qquad t>0 , \qquad x\in{\mathbb R} .
\]
We shall establish that  $F_{k}(f)$  possesses the following integral representation
\[
F_{k}(f)(t,x)=\int_{{\mathbb R}}f(y)A_{k}(t,x,y)|y|^{2k}dy ,
\]
 where
\[
A_{k}(t,x,y)=\frac{t}{\sqrt{4\pi}}\int_{0}^{+\infty}P_{k}(u,x,y)u^{-\frac{3}{2}} e^{-\frac{t^{2}}{4u}}du
\]
 is the Poisson kernel associated with $\{h_{n}^{k}\}_{n\geq0}$. We shall show that $F_{k}(f)(t, \cdot )\in L^{p}({\mathbb R},|x|^{2k}dx)$  and
\[
\|F_{k}(f)(t,\cdot)\|_{k,p}\leq 2^{k+\frac{1}{2}}e^{-t\sqrt{2k+1}}\|f\|_{k,p}.
\]
 Also, we def\/ine the Hilbert transforms associated with the Dunkl--Hermite operators formally~by
\[
{\cal H}_{k}^{\pm}=(T_{k}\pm x)(-L_{k})^{-\frac{1}{2}}.
\]
We write  $f\sim\sum\limits_{n=0}^{+\infty}a_{n}^{k}(f)h_{n}^{k}$, to say that the last series represents the expansions of $f$ with respect to  $\{h_{n}^{k}\}_{n\geq0}$. Note, that if
$f\sim\sum\limits_{n=0}^{+\infty}a_{n}^{k}(f)h_{n}^{k}$, then again formally,
\begin{gather*}
{\cal H}_{k}^{+}f\sim\sum_{n=0}^{+\infty} a_{n}^{k}(f)\dfrac{\theta(n,k)}{\sqrt{2n+2k+1}}h_{n-1}^{k} ,\qquad
{\cal H}_{k}^{-}f\sim-\sum_{n=0}^{+\infty}a_{n}^{k}(f)\frac{\theta (n+1,k)}{\sqrt{2n+2k+1}}h_{n+1}^{k},
\end{gather*}
 where
\[
\theta(n,k)=\left\{\begin{array}{ll}
\sqrt{2n} & \mbox{if} \   n \  \mbox{is even},\\
 \theta(n,k)=\sqrt{2n+4k} &\mbox{if} \ n \ \mbox{is odd};
 \end{array}\right.
\]
here and also later on, we use the convention that $h_{n-1}^{k}=0$  if  $n=0$.

We shall see that
\[
{\cal H}_{k}^{\pm}f(x)=\lim_{\epsilon\longrightarrow0}\int_{|x-y|>\epsilon} f(y)R_{k}^{\pm}(x,y)|y|^{2k}dy,
 \]
 exist for almost every $x$,
where  $R_{k}^{\pm}(x,y)$  are appropriate kernels. Next, we shall prove that the operators ${\cal H}_{k}^{\pm}$ are bounded on
$L^{p}({\mathbb R},|x|^{2k}dx)$.

 Finally, we use the Dunkl--Hermite functions to def\/ine the conjugate Poisson integrals $f_{k}^{\pm}(t,x)$  by
\begin{gather*}
f_{k}^{+}(t,x)=\sum_{n=0}^{+\infty}e^{-t\sqrt{2n+2k+1}} a_{n}^{k}(f)\dfrac{\theta(n,k)}{\sqrt{2n+2k+1}}h_{n-1}^{k}(x),\\
f_{k}^{-}(t,x)=\sum_{n=0}^{+\infty}e^{-t\sqrt{2n+2k+1}} a_{n}^{k}(f)\dfrac{\theta(n+1,k)} {\sqrt{2n+2k+1}}h_{n+1}^{k}(x).
\end{gather*}
We shall establish that $f_{k}^{\pm}(t,x)$ possesses the integral representations
\begin{gather*}
f_{k}^{+}(t,x)=\int_{{\mathbb R}}Q_{k}(t,x,y)f(y)|y|^{2k}dy,\qquad
f_{k}^{-}(t,x)=\int_{{\mathbb R}}M_{k}(t,x,y)f(y)|y|^{2k}dy,
\end{gather*}
where $Q_{k}(t,x,\cdot)$  and  $M_{k}(t,x,\cdot)$  are kernels expressed in terms of the Dunkl kernel  $E_{k}(x,y)$  which is the eigenfunction of the Dunkl operator  $T_{k}$.

We point out that recently (see \cite{Now})  A.~Nowak and K.~Stempak have studied Riesz transforms for the Dunkl harmonic oscillator in the rank-one case.

We conclude this introduction by giving the organization of this paper. In the second section, we def\/ine the heat-dif\/fusion
integral $G_{k}(f)$  and the Poisson integral  $F_{k}(f)$ and give the integral representations of $G_{k}(f)$  and $F_{k}(f)$.
In the third section, we deal with the Hilbert transforms~${\cal H}_{k}^{\pm}$  and prove that these operators are of the strong type  $(p,p)$, $1<p<+\infty$. The remaining part is concerned with the study of the conjugate Poisson.

\section[Heat-diffusion and Poisson integrals]{Heat-dif\/fusion and Poisson integrals}

\subsection[Heat-diffusion]{Heat-dif\/fusion}

With the help of the Dunkl--Hermite functions introduced in  the previous section, we def\/ine
and study the heat-dif\/fusion in the Dunkl setting.   As the Dunkl--Hermite polynomials are expressed in terms of Laguerre polynomials, using Lemma~1.5.4 in \cite{Than}, we have the following limiting behavior of $\|h_{n}^{k}\|_{k,p}$, with respect to~$n$.

\begin{proposition}\label{P1}
For  $1\leq p\leq 4$, we have
\begin{gather*}
\|h_{2n}^{k}\|_{k,p} \sim \left\{
\begin{array}{ll}
n^{-\frac{1}{4}+\frac{1}{2p}+k(\frac{1}{p}-\frac{1}{2})} , & \mbox{if} \ \ k(p-2)<1,\\
 n^{-\frac{1}{4}-\frac{1}{2p}+k(\frac{1}{2}-\frac{1}{p})} , &\mbox{if} \ \ k(p-2)>1,
 \end{array}\right.\\
\|h_{2n+1}^{k}\|_{k,p}\sim \left\{\begin{array}{ll}
 n^{-\frac{1}{4}+\frac{1}{2p}+k(\frac{1}{p}-\frac{1}{2})} , & \mbox{if} \ \ k(p-2)<1,\\
n^{-\frac{1}{4}-\frac{1}{2p}+k(\frac{1}{2}-\frac{1}{p})}, & \mbox{if} \ \ k(p-2)>1.
\end{array}\right.
\end{gather*}
For $p>4$, we have
\begin{gather*}
\|h_{2n}^{k}\|_{k,p}\sim \left\{
\begin{array}{ll}
n^{-\frac{1}{12}-\frac{1}{6p}+k(\frac{1}{p}- \frac{1}{2})}, & \displaystyle \mbox{if} \ \ k(p-2)\leq \frac{1}{3}+\frac{p}{6},\vspace{1mm}\\
 n^{-\frac{1}{4}-\frac{1}{2p}+k(\frac{1}{2}-\frac{1}{p})} , & \displaystyle \mbox{if} \ \ k(p-2)>\frac{1}{3}+\frac{p}{6},
 \end{array}\right.\\
\|h_{2n+1}^{k}\|_{k,p}\sim \left\{
\begin{array}{ll}
n^{-\frac{1}{12}-\frac{1}{6p}+k(\frac{1}{p}-\frac{1}{2})} , & \displaystyle \mbox{if} \ \ k(p-2)\leq\frac{1}{3}+\frac{p}{6},\vspace{1mm}\\
n^{-\frac{1}{4}-\frac{1}{2p}+k(\frac{1}{2}-\frac{1}{p})}, & \displaystyle \mbox{if} \ \ k(p-2)>\frac{1}{3}+\frac{p}{6}.
\end{array}\right.
\end{gather*}
\end{proposition}

\begin{proposition}\label{P2}
There exists a positive constant  $C$  such that
\[
\|h_{n}^{k}\|_{\infty}\leq C n^{\frac{k}{2}-\frac{1}{12}}.
\]
\end{proposition}

\begin{proof}  Using the fact that
\[
H_{2n}^{k}=\frac{\sqrt{n!\Gamma(n+k+\frac{1}{2})}}{(2n)!\Gamma(k+\frac{1}{2})}V_{k}(H_{2n})
 \]
 (see \cite{Rosen}), where $\{H_{2n}\}_{n\geq 0}$  is the set of classical
Hermite polynomials, and  $V_{k}$  is the intertwining operator between $T_{k}$  and the usual derivative $\frac{d}{dx}$ given by
\[
V_{k}(f)(x)=\dfrac{2^{-2k}\Gamma(2k+1)}{\Gamma(k)\Gamma(k+1)} \int_{-1}^{1}f(xt)\left(1-t^{2}\right)^{k-1}(1+t)dt,
\]
we deduce that
\[
\|h_{2n}^{k}\|_{\infty}\leq C \frac{2^{n}\sqrt{n!\Gamma(n+k+\frac{1}{2})}}{\sqrt{(2n)!}\Gamma(k+\frac{1}{2})}\|h_{2n}\|_{\infty},
\]
where  $\{h_{2n}\}_{n\geq0}$ is the classical Hermite-functions. In view of the following estimate
\[
\|h_{n}\|_{\infty}\leq Cn^{-\frac{1}{12}}
\]
 given in \cite[Lemma~2.1]{Goss}, using Stirling's formula, we can
deduce easily the result.

 In the same way, we have
\[
H_{2n+1}^{k}=\frac{\sqrt{n!\Gamma(n+k+\frac{3}{2})}}{(2n+1)!\Gamma(k+\frac{1}{2})}V_{k}(H_{2n+1}),
\]
where  $\{H_{2n+1}\}_{n\geq 0}$ is the classical Hermite polynomials, we obtain the result as above.
\end{proof}

Let
\[
a_{n}^{k}(f)=\int_{{\mathbb R}}f(t)h_{n}^{k}(t)|t|^{2k}dt.
\]
Using the fact that $h_{n}^{k}\in L^{p'}({\mathbb R},|x|^{2k}dx)$, where $p'$  is the conjugate exponent of  $p$, and H\"{o}lder's inequality, we deduce that
\[
|a_{n}^{k}(f)|\leq\|f\|_{k,p}\|h_{n}^{k}\|_{k,p'}.
\]
In view of Propositions  \ref{P1}  and  \ref{P2}, we have the following.

\begin{proposition}\label{P3}
The series
\[
\sum_{n=0}^{+\infty}e^{-t(2n+2k+1)}a_{n}^{k}(f)h_{n}^{k}(x) , \qquad t>0 ,
\]
 converges uniformly in  $x\in{\mathbb R}$.
\end{proposition}
\begin{definition}\label{e}
We def\/ine the heat-dif\/fusion integral of  $f$  by
\begin{equation*}
G_{k}(f)(t,x)=\sum_{n=0}^{+\infty}e^{-t(2n+2k+1)}a_{n}^{k}(f)h_{n}^{k}(x), \qquad t>0 , \qquad x\in{\mathbb R}.
\end{equation*}
\end{definition}

\begin{proposition}\label{P4}
The heat-diffusion integral  $G_{k}(f)$ possesses the following integral representation
\[
G_{k}(f)(t,x)=\int_{{\mathbb R}}P_{k}(t,x,y)f(y)|y|^{2k}dy,
\]
where
\[
P_{k}(t,x,y)=\sum_{n=0}^{+\infty}e^{-t(2n+2k+1)}h_{n}^{k}(x)h_{n}^{k}(y).
\]
\end{proposition}
\begin{proof} We obtain an integral form of  $G_{k}(f)$  by writing
 \begin{gather*}
 G_{k}(f)(t,x)  =\sum_{n=0}^{+\infty}e^{-t(2n+2k+1)}h_{n}^{k}(x)\int_{{\mathbb R}}f(y)h_{n}^{k}(y) |y|^{2k}dy
\\
\phantom{G_{k}(f)(t,x)}{} =\int_{{\mathbb R}}\sum_{n=0}^{+\infty}e^{-t(2n+2k+1)}h_{n}^{k}(x)h_{n}^{k}(y)f(y)|y|^{2k}dy
 =\int_{{\mathbb R}}P_{k}(t,x,y)f(y)|y|^{2k}dy.
\end{gather*}
Interchanging the order of summation and integration is justif\/ied by Lebesgue's dominated convergence theorem since
\begin{gather*}
\sum_{n=0}^{+\infty}e^{-t(2n+2k+1)}\int_{{\mathbb R}}|h_{n}^{k}(x)h_{n}^{k}(y)f(y)||y|^{2k}dy \\
\qquad{} \leq \sum_{n=0}^{+\infty}e^{-t(2n+2k+1)}\|h_{n}^{k}\|_{\infty}\|h_{n}^{k}\|_{k,p'}\|f\|_{k,p}
<+\infty.\tag*{\qed}
\end{gather*}\renewcommand{\qed}{}
\end{proof}

Mehler's formula established by Margit R\"{o}sler for Dunkl--Hermite polynomials (see \cite{Rosl0}), adapted to Dunkl--Hermite functions $\{h_{n}^{k}\}_{n\geq0}$ reads
\begin{equation}\label{(2)}
\sum_{n=0}^{+\infty}r^{n}h_{n}^{k}(y)h_{n}^{k}(z)
=\frac{c_{k}}{(1-r^{2})^{k+\frac{1}{2}}}e^{-\frac{1}{2}(\frac{1+r^{2}}{1-r^{2}})(y^{2}+z^{2})}E_{k}
\left(\frac{2ry}{1-r^{2}},z\right)  , \qquad 0<r<1,
\end{equation} where $c_{k}$ is the constant def\/ined by
\[
c_{k}=\left(\int_{{\mathbb R}}e^{-x^{2}}|x|^{2k}dx\right)^{-1}
\]
 and  $E_{k}$ is the Dunkl kernel expressed in terms of the normalized Bessel function
\[
E_{k}(x,y)=j_{k-\frac{1}{2}}(ixy)+\frac{xy}{2k+1}j_{k+\frac{1}{2}}(ixy)  ,
\]
 where
\[
j_{\alpha}(z)=\Gamma(\alpha+1)\sum_{n=0}^{+\infty}\dfrac{(-1)^{n}(\frac{z}{2})^{2n}}{n!\Gamma(n+\alpha+1)}, \qquad \alpha\geq-\frac{1}{2}.
\]
 Set
\[
U_{k}(r,y,z): =\sum_{n=0}^{+\infty}r^{n}h_{n}^{k}(y)h_{n}^{k}(z), \qquad 0<r<1.
\]

\begin{proposition}\label{P5}
The kernel  $U_{k}$  satisfies the following properties
\begin{gather}
(i)\quad U_{k}(r,y,z)\geq0,\nonumber 
\\
(ii) \quad  U_{k}(r,y,z)=U_{k}(r,z,y),
\\
(iii) \quad
 \|U_{k}(r,y,\cdot )\|_{k,1}=\left(\frac{2}{1+r^{2}}\right)^{k+\frac{1}{2}}e^{-\frac{1}{2}\big(\frac{1-r^{2}}{1+r^{2}}\big)y^{2}}.
 \label{(5)}
 \end{gather}
\end{proposition}

\begin{proof}
$(i)$  and  $(ii)$  are obvious, let us therefore prove $(iii)$.
\begin{gather*}
\|U_{k}(r,y,\cdot)\|_{k,1} = \frac{c_{k}}{(1-r^{2})^{k+\frac{1}{2}}} \int_{{\mathbb R}}e^{-\frac{1}{2}\big(\frac{1+r^{2}}{1-r^{2}}\big)(y^{2}+z^{2})}E_{k}\left(\frac{2ry}{1-r^{2}},z\right)|z|^{2k}dz
\\
\phantom{\|U_{k}(r,y,\cdot)\|_{k,1}}{} =\frac{c_{k}}{(1-r^{2})^{k+\frac{1}{2}}}e^{-\frac{1}{2}\big(\frac{1+r^{2}}{1-r^{2}}\big)y^{2}} \int_{{\mathbb R}}e^{-\frac{1}{2}\left(\frac{1+r^{2}}{1-r^{2}}\right)z^{2}}E_{k}\left(\frac{2ry}{1-r^{2}},z\right)|z|^{2k}dz.
\end{gather*}
Performing the change of variables
$
u=\sqrt{\frac{1+r^{2}}{1-r^{2}}}z,
$
 and using the following identity (see \cite{Dunkl})
\[
\int_{{\mathbb R}}E_{k}(x,y)e^{-\frac{x^{2}}{2}}|x|^{2k}dx=2^{k+\frac{1}{2}}c_{k}^{-1}e^{
\frac{y^{2}}{2}}
\]
we are done.
\end{proof}

\begin{proposition}\label{P6}
The heat-diffusion integral $G_{k}(f)$ is a $C^{\infty} $ function on ${\mathbb R}_{+}\times{\mathbb R}$ satisfying the differential-difference equation
\begin{equation}\label{(6)}
\left(L_{k,x}-\frac{\partial}{\partial t}\right)G_{k}(f)(t,x)=0,
\end{equation}
$($here  $L_{k,x}$ means that the operator  $L_{k}$ acts on the variable $x)$.
\end{proposition}
\begin{proof} On the one hand, one has for all $m\in{\mathbb N}$
\begin{equation}\label{(7)}
\dfrac{\partial^{m}}{{\partial t}^{m}}G_{k}(f)(t,x)=\sum_{n=0}^{+\infty}(-1)^{m}(2n+2k+1)^{m}e^{-t(2n+2k+1)}a_{n}^{k}(f)h_{n}^{k}(x).
 \end{equation}
On the other hand, it is easy to see that
\[
(h_{n}^{k})'(x)=e^{-\frac{x^{2}}{2}}\frac{\partial}{\partial x}H_{n}^{k}(x)-x h_{n}^{k}(x),
\]
thus for f\/ixed $t$, the series \eqref{(7)} can be dif\/ferentiated termwisely with respect to $x$. A similar argument holds for higher derivatives and then $G_{k}(f)$ is $C^{\infty}$ on
${\mathbb R}_{+}\times {\mathbb R}$. Dif\/ferentiating term by term shows that $G_{k}(f)(t,x)$ satisf\/ies \eqref{(6)}.
\end{proof}

\begin{theorem}
The heat-diffusion integral $G_{k}(f)(t,\cdot)$, $t>0$, satisfies
\[
\|G_{k}(f)(t,\cdot)\|_{k,p}\leq(\cosh(2t))^{-(k+\frac{1}{2})}\|f\|_{k,p}.
\]
\end{theorem}
\begin{proof} Using
\[
P_{k}(t,x,y)=e^{-t(2k+1)}U_{k}\big(e^{-2t},x,y\big)
\]  and  \eqref{(5)},  we obtain
\[
\int_{{\mathbb R}}P_{k}(t,x,y)|y|^{2k}dy=(\cosh(2t))^{-(k+\frac{1}{2})}e^{-\frac{1}{2}\tanh(2t)x^{2}}.
\]
By H\"{o}lder's inequality, it follows that
\[
|G_{k}(f)(t,x)|^{p}\leq(\cosh(2t))^{-\frac{p(k+\frac{1}{2})}{p'}}\int_{{\mathbb R}}|f(y)|^{p}|P_{k}(t,x,y)||y|^{2k}dy,
\]
where  $p'$  is the conjugate exponent of  $p$. Integration with respect to $x$ and using Fubini's
Theorem yield
\begin{equation*}
\|G_{k}(f)(t,\cdot)\|_{k,p}\leq(\cosh(2t))^{-\big(k+\frac{1}{2}\big)}\|f\|_{k,p}.\tag*{\qed}
\end{equation*}\renewcommand{\qed}{}
\end{proof}

\subsection{Poisson integral}
In this subsection, we introduce the Poisson integral and we give its $L^{p}$ boundedness.
\begin{definition}
The Poisson integral  $F_{k}(f)$  of  $f$  is def\/ined by
\begin{equation*}
F_{k}(f)(t,x)=\sum_{n=0}^{+\infty}e^{-t\sqrt{2n+2k+1}}a_{n}^{k}(f)h_{n}^{k}(x),\qquad t>0 , \qquad x\in{\mathbb R}.
\end{equation*}
\end{definition}
Again the def\/ining series is convergent by Propositions \ref{P1}  and  \ref{P2}.
\begin{proposition}\label{P7}
$F_{k}(f)$  possesses the following integral representation
\[
F_{k}(f)(t,x)=\int_{{\mathbb R}}f(y)A_{k}(t,x,y)|y|^{2k}dy  ,
\]
where
\begin{equation}\label{(10)}
A_{k}(t,x,y)=\frac{t}{\sqrt{4\pi}}\int_{0}^{+\infty}P_{k}(u,x,y)u^{-\frac{3}{2}}e^{-\frac{t^{2}}{4u}}du.
\end{equation}
$A_{k}$ is called the Poisson kernel associated with $\{h_{n}^{k}\}_{n\geq0}$.
\end{proposition}

\begin{proof}
By using the subordination formula
\begin{equation}\label{(11)}
e^{-\beta}=\frac{\beta}{\sqrt{4\pi}}\int_{0}^{+\infty}e^{-s}s^{-\frac{3}{2}} e^{-\frac{\beta^{2}}{4s}}ds  , \qquad \beta>0,
\end{equation}
we obtain an integral form of  $F_{k}(f)(t,x)$  by writing
\begin{gather*}
F_{k}(f)(t,x)  =\sum_{n=0}^{+\infty}\frac{t\sqrt{2n+2k+1}}{\sqrt{4\pi}} h_{n}^{k}(x)\int_{0}^{+\infty}e^{-s}s^{-\frac{3}{2}}e^{-\frac{t^{2}(2n+2k+1)}{4s}}ds
\int_{{\mathbb R}}f(y)h_{n}^{k}(y)|y|^{2k}dy
\\
\phantom{F_{k}(f)(t,x)}{}  =\frac{t}{\sqrt{4\pi}}\sum_{n=0}^{+\infty}h_{n}^{k}(x)\int_{0}^{+\infty}u^{-\frac{3}{2}}e^{-u(2n+2k+1)}e^{-\frac{t^{2}}{4u}}du\int_{{\mathbb R}}f(y)h_{n}^{k}(y)|y|^{2k}dy
\\
 \phantom{F_{k}(f)(t,x)}{}=\frac{t}{\sqrt{4\pi}}\int_{{\mathbb R}}f(y)\int_{0}^{+\infty}\sum_{n=0}^{+\infty}e^{-(2n+2k+1)u}h_{n}^{k}(x)h_{n}^{k}(y)
 u^{-\frac{3}{2}}e^{-\frac{t^{2}}{4u}}du|y|^{2k}dy
\\
 \phantom{F_{k}(f)(t,x)}{}=\frac{t}{\sqrt{4\pi}}\int_{{\mathbb R}}f(y)\int_{0}^{+\infty}P_{k}(u,x,y)u^{-\frac{3}{2}}e^{-\frac{t^{2}}{4u}}du|y|^{2k}dy
 =\int_{{\mathbb R}}f(y)A_{k}(t,x,y)|y|^{2k}dy.
\end{gather*}
The same argument used for the heat-dif\/fusion integral implies that $F_{k}(f)$  is  $C^{\infty}$ on  ${\mathbb R}_{+}\times {\mathbb R}$  and satisf\/ies
\begin{equation*}
\left(L_{k,x}+\frac{\partial^{2}}{\partial t^{2}}\right)F_{k}(f)(t,x)=0.\tag*{\qed}
\end{equation*}\renewcommand{\qed}{}
\end{proof}
\begin{theorem}
$F_{k}(f)(t,\cdot )\in L^{p}({\mathbb R},|x|^{2k}dx)$  and
\[
\|F_{k}(f)(t,\cdot )\|_{k,p}\leq 2^{k+\frac{1}{2}}e^{-t\sqrt{2k+1}}\|f\|_{k,p}.
\]
\end{theorem}

\begin{proof} One has
\begin{gather*}
F_{k}(f)(t,x)  = \int_{{\mathbb R}}A_{k}(t,x,y)f(y)|y|^{2k}dy
 = \frac{t}{\sqrt{4\pi}}\int_{{\mathbb R}}\int_{0}^{+\infty}P_{k}(u,x,y)u^{-\frac{3}{2}} e^{-\frac{t^{2}}{4u}}duf(y)|y|^{2k}dy
\\
  \phantom{F_{k}(f)(t,x)}{}
 =\frac{t}{\sqrt{4\pi}}\int_{0}^{+\infty}\left(\int_{{\mathbb R}}P_{k}(u,x,y)f(y)|y|^{2k}dy\right)u^{-\frac{3}{2}}e^{-\frac{t^{2}}{4u}}du
\\
\phantom{F_{k}(f)(t,x)}{}
 =\frac{t}{\sqrt{4\pi}}\int_{0}^{+\infty}G_{k}(f)(u,x)u^{-\frac{3}{2}}e^{-\frac{t^{2}}{4u}}du.
\end{gather*}
It follows that
\begin{equation*}
\|F_{k}(f)(t,\cdot)\|_{k,p}\leq\frac{t}{\sqrt{4\pi}}\int_{0}^{+\infty}\|G_{k}(f)(u,\cdot)\|_{k,p}
u^{-\frac{3}{2}}e^{-\frac{t^{2}}{4u}}du \leq 2^{k+\frac{1}{2}}e^{-t\sqrt{2k+1}}\|f\|_{k,p}.\tag*{\qed}
\end{equation*}
\renewcommand{\qed}{}
\end{proof}

\section{Hilbert transforms}
The operator  $(-L_{k})$  is positive and symmetric in $L^{2}({{\mathbb R}},|x|^{2k}dx)$ on the domain $C_{c}^{\infty}({{\mathbb R}})$. It may be easily checked that the operator  $(-{\cal L}_{k})$  given
by
\[
(-{\cal L}_{k})\left(\sum_{n=0}^{+\infty}a_{n}^{k}(f)h_{n}^{k}\right)=\sum_{n=0}^{+\infty}(2n+2k+1)a_{n}^{k}(f)h_{n}^{k}
\]
on the domain
\[
{\rm Dom} (-{\cal L}_{k})=\left\{f\in L^{2}({{\mathbb R}},|x|^{2k}dx): \sum_{n=0}^{+\infty}|(2n+2k+1)a_{n}^{k}(f)|^{2}<+\infty\right\}
 \]
 is a self-adjoint extension of  $(-L_{k})$, has the spectrum
$\{2n+2k+1\}$ and admits the spectral decomposition
\[
(-{\cal L}_{k})f=\sum_{n=0}^{+\infty}(2n+2k+1)a_{n}^{k}(f)h_{n}^{k}, \qquad f\in {\rm Dom}(-{\cal L}_{k}).
\]
Following \cite[p.~57]{Stein} the Hilbert transforms associated with the Dunkl--Hermite operator are formally
given by
\[
{{\cal H}_{k}}^{\pm}=(T_{k}\pm x)(-L_{k})^{-\frac{1}{2}} .
\]
Note, that if  $ f \sim \sum\limits_{n=0}^{+\infty}a_{n}^{k}(f)h_{n}^{k}$, then again formally,
\begin{gather*}
{\cal H}_{k}^{+}f\sim\sum_{n=0}^{+\infty}a_{n}^{k}(f)\frac{\theta(n,k)} {\sqrt{2n+2k+1}}h_{n-1}^{k} ,\qquad
{\cal H}_{k}^{-}f\sim-\sum_{n=0}^{+\infty}a_{n}^{k}(f)\frac{\theta(n+1,k)} {\sqrt{2n+2k+1}}h_{n+1}^{k},
\end{gather*}
where
\[
\theta(n,k)=\left\{
\begin{array}{ll}
\sqrt{2n} & \mbox{if $n$ is even},\\
\sqrt{2n+4k} & \mbox{if $n$ is odd}.
\end{array}\right.
\]
We use the convention that  $h_{n-1}^{k}=0$  if  $n=0$. It is clear that  ${\cal H}_{k}^{\pm}$  are def\/ined on $L^{2}({\mathbb R},|x|^{2k}dx)$~by
\begin{gather*}
{\cal H}_{k}^{+}f=\sum_{n=0}^{+\infty}a_{n}^{k}(f)\frac{\theta(n,k)}{\sqrt{2n+2k+1}}h_{n-1}^{k},\qquad
 {\cal H}_{k}^{-}f=-\sum_{n=0}^{+\infty} a_{n}^{k}(f)\frac{\theta(n+1,k)}{\sqrt{2n+2k+1}} h_{n+1}^{k}.
 \end{gather*}
To proceed to a deeper analysis of these def\/initions, in particular to consider ${\cal H}_{k}^{\pm}$  on a wider class of functions,
we def\/ine the kernels $R_{k}^{\pm}(x,y)$ by
\begin{gather}
R_{k}^{\pm}(x,y) =\frac{1}{\sqrt{\pi}}(T_{k,x}\pm x) \int_{0}^{+\infty}P_{k}(t,x,y)t^{-\frac{1}{2}}dt
 =\frac{1}{\sqrt{\pi}}\int_{0}^{+\infty} (T_{k,x}\pm x)P_{k}(t,x,y)t^{-\frac{1}{2}}dt.\label{(15)}
\end{gather}
It will be shown in Proposition \ref{P8}  that the second integral in \eqref{(15)} converges.

We have
\begin{gather*}
P_{k}(t,x,y) =\sum_{n=0}^{+\infty}e^{-t(2n+2k+1)}h_{n}^{k}(x)h_{n}^{k}(y)
   \\
\phantom{P_{k}(t,x,y)}{} =\dfrac{c_{k}}{2^{k+\frac{1}{2}}(\sinh(2t))^{k+\frac{1}{2}}} e^{-\frac{1}{2}\coth(2t)(x^{2}+y^{2})} E_{k}\left(\frac{x}{\sinh(2t)},y\right).
\end{gather*}
The change of variables $2t=\log(\frac{1+s}{1-s})$ furnishes a useful variant of \eqref{(15)}:
\[
R_{k}^{\pm}(x,y)=\frac{\sqrt{2}}{\sqrt{\pi}}\int_{0}^{1}(T_{k,x}\pm x)K_{s}(x,y)\left(\log\left(\frac{1+s}{1-s}\right)\right)^{-\frac{1}{2}}\frac{ds}{1-s^{2}},
\]
where
\begin{equation*}
K_{s}(x,y)=c_{k}\left(\frac{1-s^{2}}{4s}\right)^{k+\frac{1}{2}} e^{-\frac{1}{4}\big(s+\frac{1}{s}\big)(x^{2}+y^{2})}E_{k}\left(\left(\frac{1-s^{2}}{2s}\right)x,y\right)
.\end{equation*}
We then write
\[
R_{k}^{\pm}(x,y)=\frac{\sqrt{2}}{\sqrt{\pi}}[R_{k,1}(x,y)\pm R_{k,2}(x,y)],
\]
where
\[
R_{k,1}(x,y)=\int_{0}^{1}T_{k,x}K_{s}(x,y)\left(\log\left(\frac{1+s}{1-s}\right)\right)^{-\frac{1}{2}}\frac{ds}{1-s^{2}}  ,
\] and
\[
R_{k,2}(x,y)=\int_{0}^{1}x K_{s}(x,y)\left(\log\left(\frac{1+s}{1-s}\right)\right)^{-\frac{1}{2}}\frac{ds}{1-s^{2}}.
\]

\begin{proposition}\label{P8}
There exists a positive constant  $C$  such that for  $(x,y)\in \Delta^{c}=\{(x,y)\in {\mathbb R}^{2}:x\neq y\}$, the kernels
$R_{k,1}(x,y)$ and $R_{k,2}(x,y)$ satisfy
\begin{gather}\label{(17)}
|R_{k,1}(x,y)|\leq\frac{C}{|x-y|},
\\
\label{(18)}
|R_{k,2}(x,y)|\leq\frac{C}{|x-y|}.
\end{gather}
\end{proposition}
\begin{proof}
We start with proving  \eqref{(18)}. We have
\[
|R_{k,2}(x,y)|\leq C\int_{0}^{1}\beta(s)|x|e^{-\frac{1}{4}\big(s+\frac{1}{s}\big)(x^{2} +y^{2})}E_{k}\left(\left(\frac{1-s^{2}}{2s}\right)x,y\right)ds ,
\]
where we let
\[
\beta(s)=(1-s)^{k-\frac{1}{2}}s^{-\big(k+\frac{1}{2}\big)}\left(\log\left(\frac{1+s}{1-s}\right)\right)^{-\frac{1}{2}} .
\]
Using the following estimate (see~\cite{Rosl1})
\[
E_{k}\left(\left(\frac{1-s^{2}}{2s}\right)x,y\right)\leq e^{\big(\frac{1-s^{2}}{2s}\big)|xy|},
\]
the same reasoning as in  \cite[pp.~460--461]{Stem} in the classical case gives the result.
In order to estimate $R_{k,1}$, write
\[
T_{k,x}K_{s}(x,y)=-\frac{1}{2}\left[s(x+y)+\frac{1}{s}(x-y)\right]K_{s}(x,y)
\]
to see that
\[
|R_{k,1}(x,y)|\leq C\int_{0}^{1}\beta(s)\left[s|x+y|+\frac{1}{s}|x-y|\right]e^{-\frac{1}{4}\big(s+\frac{1}{s}\big)(x^{2}+y^{2})} E_{k}\left(\left(\frac{1-s^{2}}{2s}\right)x,y\right)ds,
\]
and use the same above arguments used to get the bound for $R_{k,1}$.
\end{proof}

\begin{proposition}\label{P9}
There exists a positive constant  $C$  such that for $(x,y)\in\Delta^{c}=\{(x,y)\in {\mathbb R}^{2}:x\neq y\}$, if $|x-y|\geq 2|x-x'|$, then
\begin{gather}\label{(19)}
|R_{k,1}(x,y)-R_{k,1}(x',y)|\leq\frac{C|x-x'|}{|x-y|^{2}} ,
\\
\label{(20)}
|R_{k,2}(x,y)-R_{k,2}(x',y)|\leq\frac{C|x-x'|}{|x-y|^{2}} .
\end{gather}
\end{proposition}

\begin{proof} Start with
\begin{gather*}
|R_{k,2}(x,y)-R_{k,2}(x',y)|  \leq \int_{0}^{1}\big|xK_{s}(x,y)-x'K_{s}(x',y)\big| \left(\log\left(\frac{1+s}{1-s}\right)\right)^{-\frac{1}{2}}\frac{ds}{1-s^{2}}
\\
\phantom{|R_{k,2}(x,y)-R_{k,2}(x',y)|}{}
\leq C\int_{0}^{1}\beta(s)\Bigg|xe^{-\frac{1}{4}\big(s+\frac{1}{s}\big)(x^{2}+y^{2})}E_{k}\left(\left(\frac{1-s^{2}} {2s}\right)x,y\right)
\\
\phantom{|R_{k,2}(x,y)-R_{k,2}(x',y)|\leq }{}
 -x'e^{-\frac{1}{4}\big(s+\frac{1}{s}\big)(x'^{2}+y^{2})} E_{k}\left(\left(\frac{1-s^{2}}{2s}\right)x',y\right)\Bigg|ds.
\end{gather*}
Using the following estimates (see \cite{Rosl1})
\begin{gather*}
\left|\frac{\partial}{\partial x}\left(E_{k}\left(\left(\frac{1-s^{2}}{2s}\right)x,y\right)\right)\right|\leq \left(\frac{1-s^{2}}{2s}\right)|y|E_{k}\left(\left(\frac{1-s^{2}}{2s}\right)x,y\right),\\
E_{k}\left(\left(\frac{1-s^{2}}{2s}\right)x,y\right)\leq e^{\big(\frac{1-s^{2}}{2s}\big)|xy|},
\end{gather*}
then the same reasoning as in  \cite[pp.~461--463]{Stem}  in the classical case gives the result.

Considering $R_{k,1}$,  we have
\begin{gather*}
|R_{k,1}(x,y)-R_{k,1}(x',y)|  \leq\int_{0}^{1} \big|T_{k,x}K_{s}(x,y)-T_{k,x}K_{s}(x',y)\big|\left(\log\left(\frac{1+s}{1-s}\right)\right)^{-\frac{1}{2}}
\frac{ds}{1-s^{2}}
\\
 \qquad{} \leq C\int_{0}^{1}\beta(s)\Bigg|\left[s(x+y)+\frac{1}{s}(x-y)\right]e^{-\frac{1}{4}\big(s+\frac{1}{s}\big)(x^{2} +y^{2})}E_{k}\left(\left(\frac{1-s^{2}}{2s}\right)x,y\right)
\\
\qquad\quad{} -\left[s(x'+y)+\frac{1}{s}(x'-y)\right]e^{-\frac{1}{4}\big(s+\frac{1}{s}\big) (x'^{2}+y^{2})}E_{k}\left(\left(\frac{1-s^{2}}{2s}\right)x',y\right)\Bigg|ds.
\end{gather*} The proof of \eqref{(19)} follows the same steps of the one of \eqref{(20)}.
\end{proof}

\begin{proposition}\label{P10}
There exists a positive constant $C$  such that for $(x,y)\in\Delta^{c}=\{(x,y)\in {\mathbb R}^{2}:x\neq y\}$, if  $|x-y|\geq
2|y-y'|$  then
\begin{gather}\label{(21)}
|R_{k,2}(x,y)-R_{k,2}(x,y')|\leq\frac{C|y-y'|}{|x-y|^{2}} ,
\\
\label{(22)}
|R_{k,1}(x,y)-R_{k,1}(x,y')|\leq\frac{C|y-y'|}{|x-y|^{2}}.
\end{gather}
\end{proposition}
\begin{proof} The proofs of \eqref{(21)} and  \eqref{(22)}  are quite similar to the ones of  \eqref{(19)} and \eqref{(20)}.\end{proof}

\begin{lemma}
Given $m$, $m=1,2,\dots ,$  and  $f\in C_{c}^{\infty}({\mathbb R})$  there exists  $C=C_{m,f}>0$  such that
\begin{equation*}
|\langle f,h_{n}^{k}\rangle|\leq C(2n+2k+1)^{-m},
\end{equation*}
where
\begin{equation*}
\langle f,g\rangle =\int_{{\mathbb R}}f(x)\overline{g(x)}|x|^{2k}dx.
\end{equation*}
\end{lemma}

\begin{proof}
\begin{gather*}
|\langle f,h_{n}^{k}\rangle |  = \left|\int_{\mathbb R}\!f(t)h_{n}^{k}(t)|t|^{2k}dt\right|\!  =\big|(-(2n+2k+1))^{-m}\langle L_{k}^{m}f,h_{n}^{k}\rangle\big| \leq C(2n+2k+1)^{-m}.\!\!\!\!\!\!\!\tag*{\qed}
\end{gather*}\renewcommand{\qed}{}
\end{proof}

\begin{theorem}
Let  $f,g\in C_{c}^{\infty}({\mathbb R})$  with disjoint supports. Then
\begin{equation}\label{(24)}
\langle {\cal H}_{k}^{\pm}f,g \rangle =\int_{{\mathbb R}}\int_{{\mathbb R}}R_{k}^{\pm}(x,y)f(y) \overline{g(x)}|y|^{2k}dy|x|^{2k}dx.
\end{equation}
\end{theorem}

\begin{proof} We f\/irst consider ${\cal H}_{k}^{+}$. Let
\[
f=\sum_{n=0}^{+\infty}a_{n}^{k}(f)h_{n}^{k} \qquad \mbox{and} \qquad
g=\sum_{n=0}^{+\infty}b_{n}^{k}(g)h_{n}^{k}.
\]
 Then
\[
{\cal H}_{k}^{+}f=\sum_{n=0}^{+\infty}\frac{\theta(n,k)} {\sqrt{2n+2k+1}}a_{n}^{k}(f)h_{n-1}^{k}.
 \]
 The convergence of the three series is in $L^{2}({\mathbb R},|x|^{2k}dx)$) and by Parseval's
identity
\begin{equation}\label{(25)}
\langle {\cal H}_{k}^{+}f,g\rangle =\sum_{n=0}^{+\infty}\frac{\theta(n,k)}{\sqrt{2n+2k+1}} a_{n}^{k}(f)\overline{ b_{n-1}^{k}(g)}.
\end{equation}
We will show that the right sides in \eqref{(24)} and \eqref{(25)}  coincide. Note that we can see as in Proposition~\ref{P8}  that
\[
\int_{0}^{+\infty}|(T_{k,x}+x)P_{k}(t,x,y)|t^{-\frac{1}{2}}dt\leq\frac{C}{|x-y|}.
\]
This result and the assumption made on the supports of  $f$  and $g$ show that
\begin{equation}\label{(26)}
\int_{{\mathbb R}}\int_{{\mathbb R}}\int_{0}^{+\infty} |(T_{k,x}+x)P_{k}(t,x,y)|t^{-\frac{1}{2}}dt |\overline{g(x)}f(y)||y|^{2k}dy|x|^{2k}dx<+\infty.
\end{equation}
Now,
\begin{gather*}
\int_{{\mathbb R}}\int_{{\mathbb R}}R_{k}^{+}(x,y)f(y)\overline{g(x)}
|y|^{2k}dy|x|^{2k}dx \\
=\frac{1}{\sqrt{\pi}}\int_{{\mathbb R}}\int_{{\mathbb R}}
\int_{0}^{+\infty}\{(T_{k,x}+x)P_{k}(t,x,y)\}
t^{-\frac{1}{2}}dtf(y)\overline{g(x)}|y|^{2k}dy|x|^{2k}dx\\
=\frac{1}{\sqrt{\pi}}\int_{0}^{+\infty}\!\!\!\int_{{\mathbb R}}\!\int_{{\mathbb R}}\!
\left\{(T_{k,x}+x)\!\left(\sum_{n=0}^{+\infty}e^{-t(2n+2k+1)}h_{n}^{k}(x)h_{n}^{k}(y)\!\right)\!\right\}
\overline{g(x)}f(y)
|x|^{2k}dx|y|^{2k}dyt^{-\frac{1}{2}}dt\\
=\frac{1}{\sqrt{\pi}}
\int_{0}^{+\infty}\int_{{\mathbb R}}\int_{{\mathbb R}}\left\{\sum_{n=0}^{+\infty}e^{-t(2n+2k+1)}\theta(n,k)
h_{n-1}^{k}(x)h_{n}^{k}(y)\right\}\overline{g(x)}f(y)|x|^{2k}dx|y|^{2k}dy
t^{-\frac{1}{2}}dt\\
=\frac{1}{\sqrt{\pi}}\int_{0}^{+\infty}\int_{{\mathbb R}}
\left\{\sum_{n=0}^{+\infty}e^{-t(2n+2k+1)}\theta(n,k)\overline{b_{n-1}^{k}(g)}h_{n}^{k}(y)\right\}
f(y)|y|^{2k}dyt^{-\frac{1}{2}}dt\\
=\frac{1}{\sqrt{\pi}}\int_{0}^{+\infty}\left\{\sum_{n=0}^{+\infty}
e^{-t(2n+2k+1)}
\theta(n,k)a_{n}^{k}(f)\overline{b_{n-1}^{k}(g)}\right\}t^{-\frac{1}{2}}dt
=\langle {\cal H}_{k}^{+}f,g\rangle.
\end{gather*}
Note that Fubini's theorem is justif\/ied by  \eqref{(26)}. Recall that
\begin{equation*}
{\cal H}_{k}^{-}f=-\sum_{n=0}^{+\infty}\frac{\theta(n+1,k)}{\sqrt{2n+2k+1}} a_{n}^{k}(f)h_{n+1}^{k},
\end{equation*}
then one gets similarly
\begin{gather*}
\int_{{\mathbb R}}\int_{{\mathbb R}}R_{k}^{-}(x,y)f(y)\overline{g(x)}
|y|^{2k}dy|x|^{2k}dx\\
\qquad{}=\frac{1}{\sqrt{\pi}}\int_{{\mathbb R}}\int_{{\mathbb R}}
\int_{0}^{+\infty}\{(T_{k,x}-x)P_{k}(t,x,y)\}
t^{-\frac{1}{2}}dtf(y)\overline{g(x)}|y|^{2k}dy|x|^{2k}dx \\
\qquad{} =-\sum_{n=0}^{+\infty}\frac{\theta(n+1,k)}{\sqrt{2n+2k+1}}
a_{n}^{k}(f)\overline{ b_{n+1}^{k}(g)} =\langle {\cal H}_{k}^{-}f,g\rangle.\tag*{\qed}
\end{gather*}\renewcommand{\qed}{}
\end{proof}

\begin{theorem}
For almost every $x$ in ${\mathbb R}$, the Hilbert transforms are given by
\begin{equation*}
{\cal H}_{k}^{\pm}f(x)=\lim_{\epsilon\longrightarrow0}\int_{|x-y|>\epsilon}f(y)R_{k}^{\pm}
(x,y) |y|^{2k}dy.
\end{equation*}
\end{theorem}

\begin{proof} We have
\begin{gather*}
\int_{|x-y|>\epsilon}\!\! f(y)R_{k}^{\pm}(x,y)|y|^{2k}dy  = \int_{|x-y|>\epsilon} \!\! f(y)|y|^{\frac{2k}{p}}|y|^{\frac{2k}{p'}}R_{k}^{\pm}(x,y)dy
   =\int_{|x-y|>\epsilon}\!\! g(y)W_{k}^{\pm}(x,y)dy,
\end{gather*}
where  $p'$  is the conjugate exponent of  $p$,
\begin{gather*}
g(y)=f(y)|y|^{\frac{2k}{p}}, \qquad g\in L^{p}({\mathbb R},dx),\qquad
W_{k}^{\pm}(x,y)=|y|^{\frac{2k}{p'}}R_{k}^{\pm}(x,y),
\end{gather*}
$W_{k}^{\pm}(x,y)$  are Calder\'{o}n--Zygmund kernels  (see Propositions~\ref{P8}, \ref{P9}  and  \ref{P10}).
It follows that
\[
\lim_{\epsilon\longrightarrow0}\int_{|x-y|>\epsilon}f(y)R_{k}^{\pm}(x,y)|y|^{2k}dy= \lim_{\epsilon\longrightarrow0}\int_{|x-y|>\epsilon}g(y)W_{k}^{\pm}(x,y)dy
\]
exist for almost every  $x$  (see \cite[p.~55]{Jour}).
\end{proof}

\begin{remark}
For  $f\in L^{2}({\mathbb R},|x|^{2k}dx)$,  we have
\begin{equation*}
{\cal F}_{k}({\cal H}_{k}^{\pm}f)(x) =\pm i {\cal H}_{k}^{\pm}({\cal F}_{k}f)(x),
\end{equation*} where ${\cal F}_{k}$  is the Plancherel Dunkl transform,  (see~\cite{De}).
\end{remark}

\begin{theorem}
The operators  ${\cal H}_{k}^{\pm}$  are bounded on $L^{p}({\mathbb R},|x|^{2k}dx)$, $1<p<+\infty$.
\end{theorem}

\begin{proof}
Consider the truncated operators
\begin{gather*}
{\cal H}_{\epsilon,k}^{\pm}f(x)=\int_{|x-y|>\epsilon}f(y)R_{k}^{\pm}(x,y)|y|^{2k}dy,\\
\|{\cal H}_{k}^{\pm}f\|_{k,p}^{p}  =\int_{{\mathbb R}}|{\cal H}_{k}^{\pm}f(x)|^{p}|x|^{2k}dx
  =\int_{{\mathbb R}}|\lim_{\epsilon\longrightarrow0} {\cal H}_{\epsilon,k}^{\pm}f(x)|^{p}|x|^{2k}dx
 \\
 \phantom{{\cal H}_{\epsilon,k}^{\pm}f(x)}{}   =\int_{{\mathbb R}}\lim_{\epsilon\longrightarrow 0}|{\cal H}_{\epsilon,k}^{\pm}f(x)|^{p} |x|^{2k}dx
 =\int_{R}\liminf_{\epsilon\longrightarrow0}|{\cal H}_{\epsilon,k}^{\pm}f(x)|^{p}|x|^{2k}dx
\\
\phantom{{\cal H}_{\epsilon,k}^{\pm}f(x)}{}
\leq \liminf_{\epsilon\longrightarrow0}\int_{{\mathbb R}}|{\cal H}_{\epsilon,k}^{\pm}f(x)|^{p}|x|^{2k}dx
  =\liminf_{\epsilon\longrightarrow0}\int_{{\mathbb R}}\big||x|^{\frac{2k}{p}}{\cal H}_{\epsilon,k}^{\pm}f(x)\big|^{p}dx.
\end{gather*}
Yet,
\begin{equation*}
|x|^{\frac{2k}{p}}{\cal H}_{\epsilon,k}^{\pm}f(x)=\int_{|x-y|>\epsilon} f(y)|x|^{\frac{2k}{p}}R_{k}^{\pm}(x,y)|y|^{2k}dy =\int_{|x-y|>\epsilon}g(y)Z_{k}^{\pm}(x,y)dy,
\end{equation*}
where
\begin{gather*}
g(y)=f(y)|y|^{\frac{2k}{p}}, g\in L^{p}({\mathbb R},dx),\qquad
Z_{k}^{\pm}(x,y)=|x|^{\frac{2k}{p}}|y|^{\frac{2k}{p'}}R_{k}^{\pm}(x,y).
\end{gather*}
$Z_{k}^{\pm}(x,y)$  are Calder\'{o}n--Zygmund kernels  (see Propositions  \ref{P8}, \ref{P9}  and  \ref{P10}).
Let
\[
S_{g}^{\pm}(x)=\int_{{\mathbb R}}g(y)Z_{k}^{\pm}(x,y)dy.
 \]
 The operators $S_{g}^{\pm}$  are Calder\'{o}n--Zygmund type associated with the Calder\'{o}n--Zygmund kernels  $Z_{k}^{\pm}(x,y)$ (see \cite[p.~48]{Jour}), then
\[
\sup_{\epsilon>0}\left|\int_{|x-y|>\epsilon}g(y)Z_{k}^{\pm}(x,y)dy\right|
\]
are bounded on  $L^{p}({\mathbb R},dx)$  for $1<p<+\infty$  (see \cite[p.~56]{Jour}). Consequently, there exists a positive constant $C=C_{p}$ such that if $f\in L^{p}({\mathbb R},|x|^{2k}dx)$  then
\begin{gather*}
\|{\cal H}_{k}^{\pm}f\|_{k,p}\leq C\|f\|_{k,p}.\tag*{\qed}
\end{gather*}\renewcommand{\qed}{}
\end{proof}

\begin{lemma}
There exists a positive constant  $C$ such that for  $f\in L^{1}({\mathbb R},|x|^{2k}dx)$, $\lambda>0$,  we have
\[
\int_{\{x\in {\mathbb R} : \sup_{\epsilon>0}|{\cal H}_{\epsilon,k}^{\pm}f(x)|>\lambda\}}dx\leq\frac{C}{\lambda}\|f\|_{k,1}.
\]
\end{lemma}
\begin{proof}
We have
\begin{gather*}
{\cal H}_{\epsilon,k}^{\pm}f(x)  = \int_{|x-y|>\epsilon}f(y)R_{k}^{\pm}(x,y)|y|^{2k}dy
  =\int_{|x-y|>\epsilon}g(y)W_{k}^{\pm}(x,y)dy,
\end{gather*}
where
\[
g(y)=f(y)|y|^{\frac{2k}{p}} \qquad \mbox{and}\qquad W_{k}^{\pm}(x,y)=|y|^{\frac{2k}{p'}}R_{k}^{\pm}(x,y).
\]
Let
\[
S_{g}^{\pm}(x)=\int_{{\mathbb R}}g(y)W_{k}^{\pm}(x,y)dy.
\]
The operators  $S_{g}^{\pm}$  are Calder\'{o}n--Zygmund operators associated with the Calder\'{o}n--Zygmund kernels $W_{k}^{\pm}(x,y)$  then there exists a positive constant  $C$ such that for $\lambda>0$ and  $f\in L^{1}({\mathbb R},|x|^{2k}dx)$,  we have
\begin{gather*}
\int_{\{x\in {\mathbb R}:\sup_{\epsilon>0}|{\cal H}_{\epsilon,k}^{\pm}f(x)|>\lambda\}}dx \leq\frac{C}{\lambda}\|f\|_{k,1}.\tag*{\qed}
\end{gather*}\renewcommand{\qed}{}
\end{proof}

As a by-product, we have the following.
\begin{theorem}
There exists a positive constant  $C$  such that for  $f \in L^{1}({\mathbb R},|x|^{2k}dx)$,  we have
\[
\sup_{\lambda>0}\left(\lambda\int_{\{x\in {\mathbb R}:|{\cal H}_{k}^{\pm}f(x)|>\lambda\}}dx\right)\leq C\|f\|_{k,1}.
\]
\end{theorem}

\section{Conjugate Poisson integrals}

Dunkl--Hermite functions allow  to def\/ine the conjugate Poisson integrals.
\begin{definition}
The conjugate Poisson integrals $f_{k}^{\pm}(t,x)$ of $f$ are def\/ined by
\begin{gather*}
f_{k}^{+}(t,x)=\sum_{n=0}^{+\infty}e^{-t\sqrt{2n+2k+1}}a_{n}^{k}(f) \frac{\theta(n,k)}{\sqrt{2n+2k+1}}h_{n-1}^{k}(x),
\\
f_{k}^{-}(t,x)=\sum_{n=0}^{+\infty}e^{-t\sqrt{2n+2k+1}}a_{n}^{k}(f) \frac{\theta(n+1,k)}{\sqrt{2n+2k+1}}h_{n+1}^{k}(x).
\end{gather*}
\end{definition}

\begin{remark}
The same arguments used for the heat-dif\/fusion integral show that $f_{k}^{\pm}(t,x)\in C^{\infty}({\mathbb R}_{+}\times {\mathbb R})$  and satisfy the dif\/ferential-dif\/ference equations
\begin{gather}\label{(36)}
(i)\quad \left(L_{k,x}+\frac{\partial^{2}}{{\partial t}^{2}}\right)f_{k}^{\pm} (t,x)=\pm 2 f_{k}^{\pm}(t,x),
\\
\label{(37)}
(ii) \quad (T_{k,x}\pm x)F_{k}(f)(t,x)=\mp\frac{\partial}{\partial t}f_{k}^{\pm}(t,x),
\end{gather}
where $F_{k}(f)(t,x)$ is the Poisson integral of  $f$. \end{remark} We now use  \eqref{(37)}  to f\/ind an integral formula for $f_{k}^{\pm}(t,x)$. Using the subordination
formula~\eqref{(11)}, taking $\beta=t\sqrt{2n+2k+1}$, making the change of variables $s\longrightarrow (2n+2k+1)u$, and then substituting $r=e^{-2u}$ leads to the formula
\begin{equation*}
e^{-t\sqrt{2n+2k+1}}=\int_{0}^{1}L(t,r)r^{n+k+\frac{1}{2}}dr,
\end{equation*}
where
\begin{equation*}
L(t,r)=\frac{t e^{\frac{t^{2}}{2\log r}}}{(2\pi)^{\frac{1}{2}}r(-\log r)^{\frac{3}{2}}}.
\end{equation*}
Then if  $A_{k}(t,x,y)$  denotes the Poisson kernel~\eqref{(10)}, we have
\begin{gather*}
A_{k}(t,x,y)  = \sum_{n=0}^{+\infty}e^{-t\sqrt{2n+2k+1}}h_{n}^{k}(x)h_{n}^{k}(y)
  =\sum_{n=0}^{+\infty}h_{n}^{k}(x)h_{n}^{k}(y)\int_{0}^{1}L(t,r) r^{n+k+\frac{1}{2}}dr
\\
\phantom{A_{k}(t,x,y)}{}  =\int_{0}^{1}\sum_{n=0}^{+\infty}r^{n}h_{n}^{k}(x)h_{n}^{k}(y)L(t,r)r^{k+\frac{1}{2}}dr
   =\int_{0}^{1}L(t,r)U_{k}(r,x,y)r^{k+\frac{1}{2}}dr.
\end{gather*}
Combining this and  \eqref{(2)}  we obtain
\begin{gather}
(T_{k,x}+x)A_{k}(t,x,y)\nonumber\\
\qquad{}=\frac{\sqrt{2}}{\sqrt{\pi}}c_{k}e^{-\frac{1}{2}(x^{2}+y^{2})}
\int_{0}^{1}\frac{(y-rx)te^{\frac{t^{2}}{2\log
r}}e^{-\frac{r^{2}x^{2}+r^{2}y^{2}}{1-r^{2}}}}{(-\log
r)^{\frac{3}{2}}(1-r^{2})^{k+\frac{3}{2}}}E_{k}(\frac{2rx}{1-r^{2}},y)r^{k+\frac{1}{2}}dr. \label{(39)}
\end{gather}
Now
\begin{equation}\label{(40)}
(T_{k,x}+x)F_{k}(f)(t,x)=\int_{{\mathbb R}}(T_{k,x}+x)A_{k}(t,x,y)f(y)|y|^{2k}dy.
\end{equation}
From Propositions~\ref{P1}  and  \ref{P2}, it is easy to check that $f_{k}^{+}(t,x)\longrightarrow0$ as $t\longrightarrow+\infty$ and so
\[
f_{k}^{+}(t,x)=-\int_{t}^{+\infty}\frac{\partial}{\partial t}f_{k}^{+}(u,x)du.
\]
Using  \eqref{(40)}, \eqref{(39)}  and \eqref{(37)} we f\/ind after integration
\[
f_{k}^{+}(t,x)=\int_{{\mathbb R}}Q_{k}(t,x,y)f(y)|y|^{2k}dy,
\]
where
\begin{equation*}
Q_{k}(t,x,y)=e^{-\frac{1}{2}(x^{2}+y^{2})}Q_{1,k}(t,x,y)
\end{equation*}
and
\begin{equation*}
Q_{1,k}(t,x,y)=\int_{0}^{1}\frac{(y-rx)}{(1-r^{2})^{k+2}}
e^{-\frac{r^{2}x^{2}+r^{2}y^{2}}{1-r^{2}}}E_{k}\left(\frac{2rx}{1-r^{2}},y\right)W_{1,k}(t,r)dr
\end{equation*}
 with
\[
W_{1,k}(t,r)=\frac{\sqrt{2}}{\sqrt{\pi}}c_{k}\left(\frac{1-r^{2}}{-\log r}\right)^{\frac{1}{2}}e^{\frac{t^{2}}{2\log r}}r^{k+\frac{1}{2}}.
\]
We now use  \eqref{(37)}  to f\/ind an integral formula for $f_{k}^{-}(t,x)$
\begin{gather}
 (T_{k,x}-x)A_{k}(t,x,y)\nonumber\\
 \qquad{} =\frac{\sqrt{2}}{\sqrt{\pi}}c_{k} e^{-\frac{1}{2}(x^{2}+y^{2})}\int_{0}^{1}\frac{(ry-x)te^{\frac{t^{2}}{2\log
r}}e^{-\frac{r^{2}x^{2}+r^{2}y^{2}}{1-r^{2}}}}{(-\log r)^{\frac{3}{2}}(1-r^{2})^{k+\frac{3}{2}}}E_{k}\left(\frac{2rx}{1-r^{2}},y\right)r^{k-\frac{1}{2}}dr. \label{(43)}
\end{gather}
Now
\begin{equation}\label{(44)}
(T_{k,x}-x)F_{k}(f)(t,x)=\int_{{\mathbb R}}(T_{k,x}-x)A_{k}(t,x,y)f(y)|y|^{2k}dy.
\end{equation}
The same reasoning as above gives $f_{k}^{-}(t,x)\longrightarrow0$ as $t\longrightarrow+\infty$ and so
\[
f_{k}^{-}(t,x)=-\int_{t}^{+\infty}\frac{\partial}{\partial t} f_{k}^{-}(u,x)du.
\]
Using  \eqref{(44)}, \eqref{(43)}  and \eqref{(37)} we f\/ind after integration
\[
f_{k}^{-}(t,x)=\int_{{\mathbb R}}M_{k}(t,x,y)f(y)|y|^{2k}dy ,
\]
where
\begin{equation*}
M_{k}(t,x,y)=e^{-\frac{1}{2}(x^{2}+y^{2})}M_{1,k}(t,x,y)
\end{equation*}
and
\[
M_{1,k}(t,x,y)=\int_{0}^{1}\frac{(x-ry)e^{-\frac{r^{2}x^{2}+r^{2}y^{2}}{1-r^{2}}}}{(1-r^{2})^{k+2}}
E_{k}\left(\frac{2rx}{1-r^{2}},y\right)Y_{k}(t,r)dr,
\]
with
\[
Y_{k}(t,r)=\frac{\sqrt{2}}{\sqrt{\pi}}c_{k}\left(\frac{1-r^{2}}{-\log r}\right)^{\frac{1}{2}}e^{\frac{t^{2}}{2\log r}}r^{k-\frac{1}{2}}.
\]

\begin{theorem}
There exists a positive constant  $C$  such that for $1<p<+\infty$, $f\in L^{p}({\mathbb R},|x|^{2k}dx)$,  we have
\[
\|f_{k}^{\pm}(t,\cdot )\|_{k,p}\leq C e^{-t\sqrt{2k+1}}\|f\|_{k,p} .
\]
 \end{theorem}

\begin{proof}
We have
\begin{gather*}
\|f_{k}^{\pm}(t,\cdot )\|_{k,p}=\|\pm {\cal H}_{k}^{\pm}F_{k}(f)(t,\cdot )\|_{k,p} \leq C\|F_{k}(f)(t,\cdot )\|_{k,p} \leq Ce^{-t\sqrt{2k+1}}\|f\|_{k,p}.\tag*{\qed}
\end{gather*}\renewcommand{\qed}{}
\end{proof}

\begin{theorem}
There exists a positive constant  $C$  such that for  $f\in L^{1}({\mathbb R},|x|^{2k}dx)$, we have
\[
\sup_{\lambda>0}\left(\lambda \int_{\{x\in {\mathbb R}: |f_{k}^{\pm}(t,x)|>\lambda\}}dx\right) \leq C e^{-t\sqrt{2k+1}}\|f\|_{k,1}.
\]
\end{theorem}
\begin{proof}  We have
\begin{gather*}
\sup_{\lambda>0}\left(\lambda\int_{\{x\in {\mathbb R}:
|f_{k}^{\pm}(t,x)|>\lambda\}}dx\right)=\sup_{\lambda>0}\left(\lambda\int
_{\{x\in {\mathbb R}: |\pm {\cal H}_{k}^{\pm}F_{k}(f)(t,x)|>\lambda\}}dx\right) \\
\hphantom{\sup_{\lambda>0}\left(\lambda\int_{\{x\in {\mathbb R}:
|f_{k}^{\pm}(t,x)|>\lambda\}}dx\right)} {}
\leq C\|F_{k}(f)(t,\cdot )\|_{k,1} \leq C e^{-t\sqrt{2k+1}}\|f\|_{k,1}.\tag*{\qed}
\end{gather*}\renewcommand{\qed}{}
\end{proof}

\subsection*{Acknowledgments}
The authors thank the anonymous referees for their careful reading of the manuscript and their valuable suggestions to improve the style of this paper.

\pdfbookmark[1]{References}{ref}
\LastPageEnding

\end{document}